\documentclass[twoside,11pt]{amsart}
\usepackage{latexsym}
\usepackage{graphics, epsfig}
\usepackage{amsmath, amscd, amssymb}
\newtheorem{theorem}{Theorem}[section]

\begin{document}

\title{Exotic smooth structures on
       $(2n+2l-1)\mathbb{CP}^2 \sharp (2n+4l-1)\overline{\mathbb{CP}}^2$}

\author{Jongil Park}
\address{Department of Mathematical Sciences, Seoul National University\\
         San 56-1 Sillim-dong, Gwanak-gu, Seoul 151-747, Korea}
\email{jipark@math.snu.ac.kr}

\author{Ki-Heon Yun}
\address{Department of Mathematics, Konkuk University,  Seoul 143-701,
         Korea}
\email{kyun@member.ams.org}%

\thanks{Jongil Park was supported by grant No. R01-2005-000-10625-0
        from the KOSEF. He also holds a joint appointment in the Research
        Institute of Mathematics, Seoul National University.}

\date{January 30, 2007; revised at March 27, 2007}



\begin{abstract}
 As an application of `reverse engineering' technique introduced by
 R. Fintushel, D. Park and R. Stern~\cite{FPS},
 we construct an infinite family of fake
 $(2n+2l-1)\mathbb{CP}^2 \sharp (2n+4l-1) \overline{\mathbb{CP}}^2$'s
 for all $n \ge 0,\, l \ge 1$.
\end{abstract}

\maketitle

\section{Introduction}

\markboth{Jongil Park and Ki-Heon Yun}{Exotic smooth structures on
          $(2n+2l-1)\mathbb{CP}^2 \sharp (2n+4l-1)\overline{\mathbb{CP}}^2$}

 Recently, R. Fintushel, D. Park and R. Stern introduced a new surgery,
 called `reverse engineering', and they constructed a family of
 homology $(2n-1)(S^2\times S^2)$'s by performing Luttinger
 surgeries on $\Sigma_2 \times \Sigma_{n+1}$ for any $n\ge 1$~\cite{FPS}.
 Here $(2n-1)(S^2\times S^2)$ means the connected sum
 of $2n-1$ copies of a $4$-manifold $S^2\times S^2$ and
 $\Sigma_{g}$ means a Riemann surface of genus $g$.

 In this short note,
 we present an easy way to produce an infinite family of fake
 $(2n+2l-1)\mathbb{CP}^2 \sharp (2n+4l-1) \overline{\mathbb{CP}}^2$'s
 for all $n \ge 0, \, l \ge 1$.
 The main idea is very simple: For each $n \ge 0$ and $l \geq 1$,
 we first take a symplectic fiber sum with $l$ copies of $Sym^2(\Sigma_3)$
 and $n$ copies of $\Sigma_2 \times \Sigma_2$
 along an essential Lagrangian torus, where $Sym^2(\Sigma_3)$
 is the $2$-fold symmetric product of genus $3$ Riemann surface $\Sigma_3$.
 And then we perform $5l+7n$ times $\pm 1$-Luttinger surgeries and
 one more $k$-surgery on the fiber sum $4$-manifold
 $l(Sym^2(\Sigma_3))\sharp n(\Sigma_2 \times \Sigma_2)$.
 Next we show that the resultant $4$-manifolds are simply connected,
 so that they are all homeomorphic to each other.
 Finally, by applying the same technique in~\cite{FPS} for the
 computation of Seiberg-Witten invariants,
 we get the following main result.

\begin{theorem}
\label{thm-main1}
 For each $n \ge 0$ and $l \ge 1$, there are infinitely many exotic
 smooth structures on the connected sum $4$-manifold
 $(2n+2l-1)\mathbb{CP}^2 \sharp (2n+4l-1) \overline{\mathbb{CP}}^2$'s.
\end{theorem}

{\em Remark}. As a special case of Theorem~\ref{thm-main1} above,
 we have a new family of fake
 $3\mathbb{CP}^2\sharp 5\overline{\mathbb{CP}}^2$'s and
 $3\mathbb{CP}^2\sharp 7\overline{\mathbb{CP}}^2$'s. It is unclear
 whether these are diffeomorphic to the fakes in ~\cite{A, AP, ABP, BK}.

\medskip

 {\it Acknowledgements}. The authors would like to thank Ron Fintushel and
 Ron Stern for providing their preprint~\cite{FPS} and for helpful conversation.
 Without their preprint~\cite{FPS}, this article would not be written.
 The authors also thank Doug Park for pointing out some errors
 in the previous version of this article.
 The first author also wishes to thank Banff International Research Station
 for warm hospitality during his stay at the institute.

\bigskip

\section{The main construction}
\label{sec-2}

 We first briefly review a Luttinger surgery and R. Fintushel, D. Park and
 R. Stern's construction of homology $S^2\times S^2$'s
 (\cite{ADK, BK, FPS} for details).

 Let $T$ be a torus in $X$ such that $\left[ T \right]^2 =0$ and
 $\gamma \subset T$ be a simple closed curve.
 By \emph{ $p/q$-surgery on $T$ along $\gamma$}, denoted by $(T,\gamma, p/q)$,
 we mean
 \[ X_{T, \gamma}(p/q) =
 ( X - \nu(T) ) \cup_{\varphi} (S^1 \times S^1 \times D^2).\]
 Here
 $\varphi : S^1 \times S^1 \times \partial D^2 \to \partial(X - \nu(T))$
 denotes a gluing map satisfying $\varphi([\partial D^2]) = q[\gamma'] + p[\mu_T]$
 in $H_1(\partial(X - \nu(T))$, where $\mu_T$ is a meridian of $T$ and
 $\gamma' \subset \partial \nu(T)$ is a simple closed curve homologous to $\gamma$
 in a tubular neighborhood $\nu(T)$ of $T$.
 If $T$ is an embedded Lagrangian torus in a symplectic
 4-manifold $X$ and $\gamma \subset T$ is a co-oriented simple closed
 curve, then, by the Darboux-Weinstein theorem,
 each nearby torus $T'$ of $T$  in $\nu(T)$ is also Lagrangian torus and we can
 consider $\gamma'$ as Lagrangian push off of $\gamma$.
 In particular, when $p=1$ and $q=k$, the $1/k$-surgery on $T$ along $\gamma$
 is called \emph{$1/k$-Luttinger surgery on $T$ along $\gamma$}.
 Some of well-known properties of $1/k$-Luttinger surgery on $T$ along
 $\gamma$ are the following:
\begin{enumerate}
 \item $\pi_1(X_{T,\gamma}(p/q)) = \pi_1(X- T)/N(\mu_T^p \gamma'^q)$.
 \item $\sigma(X)=\sigma(X_{T,\gamma}(p/q))$ and $e(X)=e(X_{T,\gamma}(p/q))$.
 \item If $X$ is a symplectic manifold, then so is $X_{T, \gamma}(1/k)$.
 \item If $\gamma^*$ is the loop $\gamma$ with the opposite co-orientation,
       then $X_{T,\gamma^*}(1/k)$ is symplectomorphic to
       $X_{T,\gamma}((-1)/k)$.
\end{enumerate}

\noindent

 By using these properties above, R. Fintushel, D. Park and R. Stern
 constructed homology $S^2\times S^2$'s by performing eight times
 $\pm 1$-Luttinger surgeries on disjoint Lagrangian tori in
 $\Sigma_2 \times \Sigma_2$:
 \begin{displaymath}
 \begin{array}{llll}
 (a'_1 \times c'_1, a'_1, -1), & (b'_1 \times c''_1, b'_1, -1),&
 (a'_2 \times c'_2, a'_2, -1), & (b'_2 \times c''_2, b'_2, -1), \\
 (a'_2 \times c'_1, c'_1, +1), & (a''_2 \times d'_1, d'_1, +1), &
 (a'_1 \times c'_2, c'_2, +1), &  (a''_1 \times d'_2, d'_2, +k)
 \end{array}
 \end{displaymath}
 where $\{ a_{1}, b_{1}, a_{2}, b_{2} \}$, $\{c_{1}, d_{1}, c_{2}, d_{2}\}$
 are the standard homotopy generators of $\Sigma_2$ based at $x$ and $y$
 respectively and $a_i'$, $a_i''$ are parallel copies of $a_i$ in $\Sigma_2$
 such that, when we consider them as  based loops, $a_i'$ is homotopic to $a_i$
 and $a_i''$ is homotopic to $b_i a_i b_i^{-1}$ relative to the base point
 $\{x\}\times\{y\}$.
 (Similarly, we choose parallel copies $b_i'$, $b_i''$, $b_i'$, $b_i''$,
 $d_i'$, $d_i''$ for $i=1, 2$.)
 Then there are $18$ relations (\cite{FPS} for details):
 \begin{align*}
  &[b_1^{-1}, d_1^{-1}] = a_1,
  &&[a_1^{-1}, d_1] = b_1,
  &&[b_2^{-1}, d_2^{-1}] = a_2,
  &&[a_2^{-1}, d_2] = b_2, \\
  &[d_1^{-1}, b_2^{-1}] = c_1,
  &&[c_1^{-1}, b_2] = d_1,
  &&[d_2^{-1}, b_1^{-1}] = c_2,
  &&[c_2^{-1}, b_1] = d_2,\\
  &[a_1, c_1] =1,
  &&[a_1, c_2] = 1,
  &&[a_1, d_2] = 1,
  &&[b_1, c_1] =1, \\
  &[a_2, c_1] =1,
  &&[ a_2, c_2] =1,
  &&[a_2, d_1] = 1,
  &&[b_2, c_2] =1, \\
  &[a_1, b_1][a_2, b_2] = 1,
  &&[c_1, d_1][c_2,d_2] =1 &&&
 \end{align*}
 in the fundamental group and it gives a family of homology $S^2 \times S^2$.

\medskip

\begin{proof}[Proof of Theorem~\ref{thm-main1}]
 We first construct an infinite family of homology
 $(2n+2l-1)\mathbb{CP}^2 \sharp (2n+4l-1)\overline{\mathbb{CP}}^2$'s
 by using a symplectic fiber sum and a Luttinger surgery. And then we will show
 that they are in fact all simply connected.

 {\em Construction} -
 Let us consider $n$ copies of
 $\Sigma_2 \times \Sigma_2$ with fixed base points
 $(x_i, y_i)$ such that $x_1=x_2=\cdots =x_n \in \Sigma_2$ and
 $y_1=y_2=\cdots =y_n \in \Sigma_2$,
 and let $\{a_{i,j}, a'_{i,j}, a''_{i,j},\cdots, d_{i,j},
 d'_{i,j}, d''_{i,j}\}$ be the simple closed curves
 in the $i$-th copy of $\Sigma_2 \times \Sigma_2$ such that $a_{i,j}'$,
 $a_{i,j}''$ are parallel copies of $a_{i,j}$ and  $a_{i,j}'$ is homotopic to
 $a_{i,j}$ and $a_{i,j}''$ is homotopic to $b_{i,j} a_{i,j} b_{i,j}^{-1}$
 relative to the base point $(x_i,y_i)$.
 (Similarly, we choose parallel copies $b_{i,j}'$, $b_{i,j}''$, $b_{i,j}'$,
 $b_{i,j}''$, $d_{i,j}'$, $d_{i,j}''$ for $j=1, 2$.)

 We  define a symplectic $4$-manifold $Z_n$ inductively as follows:
 Let us denote $Z_1 = \Sigma_2 \times \Sigma_2$.
 Note that, when performing a symplectic fiber sum, we locally perturb the
 symplectic form so that $a_{i,1}'' \times d_{i,2}'$ becomes a symplectic torus
 of square $0$ while all other disjoint Lagrangian tori are still
 Lagrangian and $a_{i,1}'' \times d_{i,2}''$ will be a symplectic torus
 of square $0$ in the $i$-th copy of $\Sigma_2 \times \Sigma_2$.
 Let $Z_i = Z_{i-1}\sharp_{a_{i-1,1}''\times d_{i-1,2}'' =  a_{i,1}'' \times
 d_{i,2}' } \Sigma_2 \times \Sigma_2$ be a symplectic $4$-manifold
 obtained by taking a symplectic fiber sum along symplectic torus
 $a_{i-1,1}'' \times d_{i-1,2}''$ in $Z_{i-1}$ and
 $a_{i,1}'' \times d_{i,2}'$ in the $i$-th copy of $\Sigma_2 \times \Sigma_2$.

 Let $Sym^2(\Sigma_3)$ be the $2$-fold symmetric product of a genus $3$ Riemann
 surface $\Sigma_3$, i.e. the quotient space of $\Sigma_3\times \Sigma_3$
 by using the
 involution $\tau : \Sigma_3 \times \Sigma_3 \to  \Sigma_3 \times \Sigma_3$
 defined by $\tau(v,w) = (w,v)$. Let $x \in \Sigma_3$ be a fixed base point
 of $\Sigma_3$ and let $f_i, g_i$, $i=1,2,3$, be standard generators of
 $\pi_1(\Sigma_3, x)$. Let $z\in Sym^2(\Sigma_3)$ be the image of $(x,x)$ and
 we will consider it as a fixed base point of $Sym^2(\Sigma_3)$.
 Then $\pi_1(Sym^2(\Sigma_3), z)= \mathbb{Z}^6$ and
 $\{f_i = f_i \times \{x\}, \, g_i = g_i \times \{x\}|\ i=1,2,3 \}$
 are generators of $\pi_1(Sym^2(\Sigma_3), z)$.

 We define $S_1 = Sym^2(\Sigma_3)$ and
 $S_s = S_{s-1} \sharp_{f_{s-1,1}'' \times g_{s-1,2}''=f_{s,1}'' \times
 g_{s,2}'} Sym^2(\Sigma_3)$  for  $s \ge 2$ by performing a symplectic fiber sum
 along torus $f_{s-1,1}'' \times g_{s-1,2}''$ and $f_{s,1}'' \times
 g_{s,2}'$ where $\{ f_{i,1}, g_{i,1}, f_{i,2}, g_{i,2}, f_{i,3}, g_{i,3} \}$
 are generators of the $s$-th copy of $\pi_1(Sym^2(\Sigma_3))$ corresponding to
 the base point $z_s$ and we consider $\{f_{s,j}', f_{s,j}'', g_{s,j}',
 g_{s,j}''\}$ as based loops in $Sym^2(\Sigma_3)$ such that $f_{s,j}'$ is
 homotopic to $f_{s,j}$ and $f_{s,j}''$ is homotopic to
 $g_{s,j} f_{s,j} g_{s,j}^{-1}$ (\cite{FPS} for details).

 For each integer $n \ge 0$ and $l \ge 1$, we define a symplectic $4$-manifold
 $Y_{l,n}$ by
 \[Y_{l,0} = S_l, \ \ Y_{l,n} = S_l \sharp_{f_{s,1}''\times
   g_{s,2}'' = a_{1,1}'' \times d_{1,2}'} Z_n .\]
 Then we construct a desired family of fake $4$-manifolds, denoted by
 $\widetilde{Y}_{l,n,k}$, from $Y_{l,n}$ by performing $7n+5l$ times
 $\pm 1$-Luttinger surgeries and one more surgery as follows:
 $7$ times $\pm 1$-Luttinger surgeries on each copy of
 $\Sigma_2 \times \Sigma_2$ -
\begin{displaymath}
 \begin{array}{llll}
 (a'_{i,1} \times c'_{i,1}, a'_{i,1}, -1), & (b'_{i,1} \times c''_{i,1},
 b'_{i,1}, -1),&
 (a'_{i,2} \times c'_{i,2}, a'_{i,2}, -1), & (b'_{i,2} \times c''_{i,2},
 b'_{i,2}, -1), \\
 (a'_{i,2} \times c'_{i,1}, c'_{i,1}, +1), & (a''_{i,2} \times d'_{i,1},
 d'_{i,1}, +1), &
 (a'_{i,1} \times c'_{i,2}, c'_{i,2}, +1), &
 \end{array}
\end{displaymath}
 $5$ times $\pm 1$-Luttinger surgeries on each copy of $Sym^2(\Sigma_3)$ -
\begin{displaymath}
 \begin{array}{lll}
  (f_{i,1}' \times f_{i,2}', f_{i,2}', +1), &
  (f_{i,1}' \times f_{i,3}'. f_{i,1}', -1), &
  (g_{i,1}' \times f_{i,3}'', g_{i,1}', -1), \\
  (f_{i,2}' \times f_{i,3}', f_{i,3}', +1), &
  (f_{i,2}'' \times g_{i,3}', g_{i,3}', +1), &
 \end{array}
\end{displaymath}
 and one more $+k/1$-surgery on $f_{1,1}''\times g_{1,2}'$ along $g_{1,2}'$.

\medskip

{\em Computation of $\pi_1(\widetilde{Y}_{l,n,k})$} - We take a base
 point $(x_i,y_i)$ in each $i$-th copy of \mbox{$\Sigma_2 \times \Sigma_2$} and
 $z_s$ for the $s$-th copy of $Sym^2(\Sigma_3)$.
 Note that $\pi_1(\Sigma_2 \times \Sigma_2, (x_i,y_i))$ is generated by
 $\{a_{i,1},\  a_{i,2},\  b_{i,1},\  b_{i,2},\  c_{i,1},\  c_{i,2},\
 d_{i,1},\ d_{i,2}\}$ and $\pi_1(Sym^2(\Sigma_3), z_s)$ is generated by
 $\{f_{s,j}, g_{s,j} | j=1,2,3\}$.

 Now we first find a group presentation of $\pi_1(\widetilde{Y}_{l,n,k}, z_1)$.
 To do this we have to find a path from $z_1$ to a base point of $\Sigma_2
 \times \Sigma_2$ or $Sym^2(\Sigma_3)$ located in $Y_{l,n}$.
 We may assume that $z_s\in \partial \nu(f_{s,1}''\times g_{s,2}')$, the fixed
 base point of $s$-th copy of $Sym^2(\Sigma_3)$.
 Since the complement of six disjoint Lagrangian tori in
 $Sym^2(\Sigma_3)$, on which we perform a surgery, is path connected,
 there are paths $\eta_{s}$, $1 \le s \le l$, from $z_s$ to a fixed point
 $Q_{s} \in \partial \nu(f_{s,1}'' \times g_{s,2}'')$.
 We may assume that $\eta_{s}$ is located on  $g_{s,1} \times f_{s,2}$ and
 a symplectic fiber sum identifies $Q_{s}$ and $z_{s+1}$, $1\le s < l$.
 Similarly, since the complement of eight disjoint Lagrangian tori in
 $\Sigma_2 \times \Sigma_2$, on which we perform a surgery, is path connected,
 there are paths $\eta_{l+i}$, $1 \le i \le n$, from $(x_i, y_i)$ to a fixed
 point $Q_{l+i} \in \partial \nu(a_{i,1}'' \times d_{i,2}'')$ such that $\eta_
 {l+i}$ is located on $b_{i,1} \times c_{i,2}$. We can also
 assume that $(x_i, y_i) \in \partial \nu(a_{i,1}''\times d_{i,2}')$ and we
 identify $Q_{l+i}$ and $(x_{i+1}, y_{i+1})$, $i=0, 1, \cdots, n-1$, when
 performing a symplectic fiber sum.

 \begin{figure}[htb]
 \includegraphics[scale = 1]{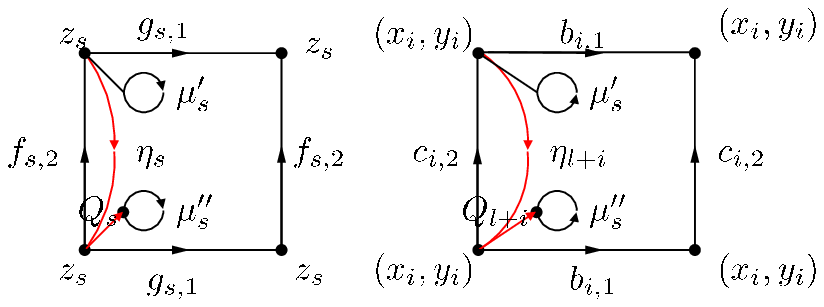}
 \caption{}\label{fig1}
 \end{figure}

 Let $\gamma_i$ be a path from $z_1$ to $z_s$ or $(x_i, y_i)$  which is
 defined as follows:
 $\gamma_1$ is a constant path based at $z_1$ and
 $\gamma_s = \eta_{1}\cdot \eta_{2} \cdots \eta_{s-1}$ for $2 \le s \le l$.
 In the same way, $\gamma_{l+i} = \eta_{1}\cdot \eta_{2} \cdots
 \eta_{l}\cdot \eta_{l+1} \cdots \eta_{l+i-1}$ be a path
 connecting $z_1$ and $(x_i,y_i)$.

 In this article, we use the notation $\alpha(\beta) = \alpha\cdot
 \beta \cdot \alpha^{-1}$ for two paths $\alpha$, $\beta$.
 Note that each symplectic fiber sum between $Sym^2(\Sigma_3)$
 gives three relations
\begin{eqnarray}
\label{eqn:rel-4}
 & & \gamma_{s-1}(g_{s-1,1}\ f_{s-1,1}\ g_{s-1,1}^{-1})
 = \gamma_{s} (g_{s,1}\ f_{s,1}\ g_{s,1}^{-1}), \\
 & & \gamma_{s-1}(f_{s-1,2}\ g_{s-1,2}\ f_{s-1,2}^{-1}) = \gamma_s(g_{s,2}),
 \nonumber \\
 & & \gamma_{s-1}(\mu_{s-1}'') =\gamma_{s}(\mu_{s}'^{-1}), \nonumber
\end{eqnarray}
 as loops based at $z_1$, where $\mu_{s-1}'$ is a meridian of
 $f''_{s-1,1}\times g'_{s-1,2}$ and $\mu_{s-1}''$ is a meridian of
 $f''_{s-1,1}\times g''_{s-1,2}$ which are considered as loops based at
 $z_{s-1}$.
 When we perform a symplectic fiber sum between $S_l$ and $\Sigma_2\times
 \Sigma_2$, we have to add the following relations:
\begin{eqnarray}
\label{eqn:rel-3}
 & & \gamma_{l+1}(b_{1,1}\ a_{1,1}\ b_{1,1}^{-1}) =  \gamma_{l}(g_{l,1}\
 f_{l,1}\ g_{l,1}^{-1}), \\
 & & \gamma_{l+1}( d_{1,2} ) = \gamma_{l}(f_{l,2}\ g_{l,2}\ f_{l,2}^{-1}), \
\ \gamma_l(\mu_{l}'') = \gamma_{l+1}(\mu_{l+1}'^{-1}). \nonumber
\end{eqnarray}

\noindent
 Similarly, performing a symplectic fiber sum between $Y_{l,i-1}\ \ ( i\ge 2)$
 and $\Sigma_2 \times \Sigma_2$, the following three relations hold:
\begin{eqnarray}
\label{eqn:rel-2}
 & & \gamma_{l+i-1}(b_{i-1,1}\  a_{i-1,1}\  b_{i-1,1}^{-1}) =
 \gamma_{l+i}( b_{i,1}\ a_{i,1}\  b_{i,1}^{-1}),  \\
 & & \gamma_{l+i-1}(c_{i-1,2}\ d_{i-1,2}\ c_{i-1,2}^{-1}) =
\gamma_{l+i}(d_{i,2}),
 \nonumber \\
 & & \gamma_{l+i-1}(\mu_{l+i-1}'') = \gamma_{l+i}(\mu_{i+i}'^{-1}),  \nonumber
\end{eqnarray}
 where $\mu_{l+i-1}'$ is a meridian of $a''_{i-1,1} \times d'_{i-1,2}$ and
$\mu_{l+i-1}''$ is a meridian of $a''_{i-1,1} \times d''_{i-1,2}$ which are
 considered as a loop based at $(x_{i-1}, y_{i-1})$.

 Note that, after performing $7$ times $\pm 1$-Luttinger surgeries
 in each copy of $\Sigma_2 \times \Sigma_2$ lying in $Y_{l,n}$,
 we have $17$ relations:
\begin{align}\label{eqn:rel-1}
  &[b_{i,1}^{-1}, d_{i,1}^{-1}] = a_{i,1},
  &&[a_{i,1}^{-1}, d_{i,1}] = b_{i,1},
  &&[b_{i,2}^{-1}, d_{i,2}^{-1}] = a_{i,2}, \\
  &[a_{i,2}^{-1}, d_{i,2}] = b_{i,2},
  &&[d_{i,1}^{-1}, b_{i,2}^{-1}] = c_{i,1},
  &&[c_{i,1}^{-1}, b_{i,2}] = d_{i,1}, \nonumber\\
  &[d_{i,2}^{-1}, b_{i,1}^{-1}] = c_{i,2},
  &&[a_{i,1}, c_{i,1}] =1,
  &&[a_{i,1}, c_{i,2}] = 1, \nonumber\\
  &[a_{i,1}, d_{i,2}] = 1,
  &&[b_{i,1}, c_{i,1}] =1,
  &&[a_{i,2}, c_{i,1}] =1, \nonumber\\
  &[ a_{i,2}, c_{i,2}] =1,
  &&[a_{i,2}, d_{i,1}] = 1,
  &&[b_{i,2}, c_{i,2}] =1,  \nonumber \\
  &[a_{i,1}, b_{i,1}][a_{i,2}, b_{i,2}] = 1,
  &&[c_{i,1}, d_{i,1}][c_{i,2},d_{i,2}] =1,  &&& \nonumber
 \end{align}
 as loops based at $(x_i, y_i)$.
 Furthermore, after performing $5$ times $\pm 1$-Luttinger surgeries in each
 copy of $Sym^2(\Sigma_3)$ lying in $Y_{l, n}$, we also have $14$ relations:
\begin{align}
  &\left[g_{s,1}^{-1}, g_{s,2}^{-1}\right] = f_{s,2},
  &&\left[g_{s,1}^{-1}, g_{s,3}^{-1}\right] = f_{s,1},
  &&\left[f_{s,1}^{-1}, g_{s,3}\right] = g_{s,1},
  &&\left[g_{s,2}^{-1}, g_{s,3}^{-1}\right] = f_{s,3}, \\
  &\left[g_{s,2}, f_{s,3}^{-1}\right] = g_{s,3},
  &&\left[f_{s,1}, g_{s,1}\right] = 1,
  &&\left[f_{s,1}, f_{s,2}\right] =1,
  &&\left[f_{s,1}, g_{s,2}\right] = 1,  \nonumber\\
  &\left[f_{s,1}, f_{s,3}\right] = 1,
  &&\left[g_{s,1}, f_{s,3}\right] = 1,
  &&\left[f_{s,2}, g_{s,2}\right] = 1,
  &&\left[f_{s,2}, f_{s,3}\right] = 1, \nonumber\\
  &\left[f_{s,2}, g_{s,3}\right] = 1,
  &&\left[f_{s,3}, g_{s,3}\right] = 1, &&& \nonumber
\end{align}
 as loops based at $z_s$ and it gives
\begin{eqnarray}
\label{eqn:rel-sym}
 & & \gamma_{s}(f_{s,1}) = \gamma_{s}(f_{s,2}) = \gamma_{s}(g_{s,1}) =1, \\
 & & \gamma_{s}(\left[g_{s,2}^{-1}, g_{s,3}^{-1}\right]) = \gamma_{s}(f_{s,3}),
\ \
 \gamma_{s}(\left[g_{s,2}, f_{s,3}^{-1}\right]) = \gamma_{s}(g_{s,3}) \nonumber
\end{eqnarray}
 as loops based at $z_1$ by the same method as in~\cite{FPS}, i.e.
 \[g_{s,1} = [f_{s,1}^{-1}, g_{s,3}] = [f_{s,1}^{-1}, [g_{s,2},
   f_{s,3}^{-1}]] =1 \]
 because $[f_{s,1}, f_{s,3}]=1$ and $[f_{s,1}, g_{s,2}]=1$.

 Note that, since $g_{1,1} =1$ due to the relation (\ref{eqn:rel-sym}),
 the last surgery $(f_{1,1}'' \times g_{1,2}', g_{1,2}', +k)$ gives a relation
\begin{equation}\label{eqn:rel-5}
 g_{1,2} =  ([g_{1,1}, f_{1,2}^{-1}])^k = 1,
\end{equation}
 so that the relations (\ref{eqn:rel-4}) and (\ref{eqn:rel-5}) give
 $\gamma_s(g_{s,2}) =1$ for each $s=2, \cdots, l$ and if we
 apply this to the relation (\ref{eqn:rel-sym}), we also get
\begin{equation}\label{eqn:sym}
 \gamma_{s}(f_{s,m}) = \gamma_{s}(g_{s,m}) =1 \textrm{ for all } s=1,
 \cdots, l,\ \  m=1, 2, 3.
\end{equation}
 Hence the relations (\ref{eqn:sym}) and (\ref{eqn:rel-3}) imply
 \[\gamma_{l+1}(a_{1,1}) = 1 = \gamma_{l+1}(d_{1,2}) \]
 as loops  based at $z_1$ and if we apply this to the relation (\ref{eqn:rel-1})
 again, then we also get
\begin{equation}\label{eqn:Z1}
 \gamma_{l+1}(a_{1,j}) = \gamma_{l+1}(b_{1,j})=\gamma_{l+1}(c_{1,j}) =
 \gamma_{l+1}(d_{1,j}) =1 \textrm{ for } j=1,2.
\end{equation}

 Now, by an induction argument using relations (\ref{eqn:Z1}), (\ref{eqn:rel-2})
 and(\ref{eqn:rel-1}),  we get
\begin{equation}\label{eqn:Z}
 \gamma_{l+i}(a_{i,j}) = \gamma_{l+i}(b_{i,j}) = \gamma_{l+i}(c_{i,j}) =
 \gamma_{l+i}(d_{i,j}) =1 \textrm{ for all } i=1, \cdots, n,\ j=1,2.
\end{equation}

 Therefore we have
 $\pi_1(\widetilde{Y}_{l,n,k}, z_1) = 1$ from the relations (\ref{eqn:sym})
 and (\ref{eqn:Z}).

\medskip

{\em Exotic smooth structures} -
 Since $e(Sym^2(\Sigma_3)) = 6$ and $\sigma(Sym^2(\Sigma_3)) = -2$,
 we can get easily
 $e(Y_{l,n}) = 4n + 6l, \ \sigma(Y_{l,n}) = -2l $ and
 $ e(\widetilde{Y}_{l,n,k}) = 4n + 6l, \ \sigma(\widetilde{Y}_{l,n,k}) = -2l,
   \ b_1(\widetilde{Y}_{l,n,k}) =0$, so that we have
 $b_2^+(\widetilde{Y}_{l,n,k})) = 2n+2l-1$,
 and $b_2^-(\widetilde{Y}_{l,n,k})) = 2n+4l-1$.
 Furthermore, since there are at least $3l$ disjoint tori of square $-1$ in
 $\widetilde{Y}_{l,n,k}$, descended from each copy of
 $Sym^2(\Sigma_3)$ (refer to Section 3 in~\cite{FPS}),
 the manifold  $\widetilde{Y}_{l,n,k}$ is nonspin,
 so that it is homeomorphic to the connected sum
 $(2n+2l-1) \mathbb{CP}^2\sharp (2n+4l-1)\overline{\mathbb{CP}}^2$.

 Finally, by applying R. Fintushel, D. Park and R. Stern's argument
 (refer to page 9 and Theorem 1 in~\cite{FPS}) to $\widetilde{Y}_{l,n,1}$,
 we conclude that a family of $4$-manifolds
 $\{\widetilde{Y}_{l,n,k}\ | \ k \ge 1 \}$
 contain infinitely many pairwise non-diffeomorphic fake
 $(2n+2l-1) \mathbb{CP}^2\sharp (2n+4l-1)\overline{\mathbb{CP}}^2$'s.
\end{proof}

{\em Remark}. Note that the fundamental group
 $\pi_1(\widetilde{Y}_{l,n,k})$ is trivial due to the presence of
 $Sym^2(\Sigma_3)$. But, in the case $l=0$, we do not know whether
 it is trivial or not as mentioned in~\cite{FPS}.


\bigskip
\bigskip


\begin{thebibliography}{BDF}

\bibitem[A]{A} A. Akhmedov, {\it Small exotic 4-manifolds}, math.GT/0612130
\bibitem[ABP]{ABP} A. Akhmedov, R. Baykur and D. Park, {\it Constructing
        infinitely many smooth structures on small $4$-manifolds},
         math.GT/0703480
\bibitem[ADK]{ADK} D. Auroux and S.K. Donaldson and L. Katzarkov,
        {\it Luttinger surgery along Lagrangian tori and non-isotopy
        for singular symplectic plane curves}, Math. Ann. 326 (2003), 185-203
\bibitem[AP]{AP} A. Akhmedov and D. Park, {\t Exotic smooth structures
        on small $4$-manifolds}, math.GT/0701664
\bibitem[BK]{BK} S. Baldridge and P. Kirk, {\t Constructions of small
        symplectic 4-manifolds using Luttinger surgery}, math.GT/0703065
\bibitem[FPS]{FPS} R. Fintushel, D. Park and R. J. Stern,
        {\it Reverse engineering small 4-manifolds}, math.GT/0701846
\end{thebibliography}
\end{document}